\documentclass{amsart}[12pt]
\usepackage{amsmath,amssymb,amscd}
\input xy
\xyoption{all}
\DeclareMathOperator{\tor}{Tor}
\DeclareMathOperator{\im}{im}
\DeclareMathOperator{\supp}{Supp}
\DeclareMathOperator{\reg}{reg}

\newtheorem*{mainthm}{Theorem}
\newtheorem{thm}{Theorem}[section]
\newtheorem{cor}[thm]{Corollary}
\newtheorem{lem}[thm]{Lemma}
\newtheorem{prop}[thm]{Proposition}
\newtheorem{defin}[thm]{Definition}

\newcommand{\p}[1]{\ensuremath{\mathbb{P}^#1}}
\newcommand{\sheaf}[1]{\ensuremath{\mathcal{#1} }}
\newcommand{\struct}[1]{\ensuremath{\sheaf{O}_{#1}}}

\newcommand{\ses}[3]{\ensuremath{0 \longrightarrow #1 \longrightarrow #2 \longrightarrow #3 \longrightarrow 0} }
\newcommand{\coh}[3]{\ensuremath{\mathrm{H}^{#1}(#2, #3) } }

\newcommand{\rmk}{\noindent{\bf Remark:  }}
\newcommand{\cok}{\frac{\sheaf{I}_L \cdot \sheaf{J}_L'}{\sheaf{J}}}
\renewcommand{\to}{\longrightarrow}
\begin{document}

\title{On the Castelnuovo-Mumford regularity of products of ideal sheaves}
\author{Jessica Sidman}
\address{Jessica Sidman, University of Michigan, Ann Arbor, MI 48109-1109}
\email{jsidman@umich.edu}
\begin{abstract}
In this paper we give bounds on the Castelnuovo-Mumford regularity of products of ideals and ideal sheaves.  In particular, we show that the regularity of a product of ideals is bounded by the sum of the regularities of its factors if the corresponding schemes intersect in a finite set of points.  We also show how approximations of sheaves can be used to bound the regularity of an arrangement of two-planes in projective space.
\end{abstract}
     
\maketitle

\section*{Introduction}
Let $S = k[x_0, \ldots, x_n]$, where $k$ is an arbitrary field.

The Castelnuovo-Mumford regularity of a module $M$ is an algebro-geometric notion that can be interpreted as giving a guide to the size of computations involving $M$.  In practice, computational methods are primarily restricted to manipulating rings, ideals, and modules, and one wants to compute the regularity of these algebraic objects.  However, one may also work geometrically with the corresponding definition of regularity for sheaves.  Somewhat surprisingly, sheaf-theoretic techniques yield statements about the regularity of products of ideals in certain situations:
\begin{mainthm}[1.8]
If $I$ and $J$ are homogeneous ideals of $S$ defining schemes in \p{n} whose intersection is a finite set of points, then $\reg(I \cdot J) \leq \reg(I)+\reg(J)$.
\end{mainthm}
\noindent Similar techniques can be used to prove that the ideal of an arrangement of $d$ two-planes in \p{n} is $d$-regular (Theorem \ref{thm: one-intersect}), which is a special case of a question of Sturmfels asking if $d$-regularity should hold regardless of dimension.

Interest in the regularity of powers of ideals was sparked by a bound for powers of ideal sheaves in \cite{bertram-ein-lazarsfeld}.  In recent years, many people have studied the regularity of powers and products of ideals from the algebraic point of view. From \cite{chandler} and \cite{geramita-gimigliano-pitteloud} we know that if $\dim(S/I) \leq 1$ then $\reg(I^m) \leq m \cdot \reg(I)$  and there are examples showing that this bound can fail if hypotheses on the dimension of $S/I$ are removed entirely \cite{sturmfels}, \cite{conca}.  Good bounds do exist for powers of monomial ideals \cite{smith-swanson} and asymptotically for powers of arbitrary ideals \cite{cutkosky-herzog-trung}, \cite{kodiyalam}, \cite{cutkosky-ein-lazarsfeld}.  More generally, Conca and Herzog consider the regularity of products of ideals and of ideals with modules in \cite{conca-herzog}.  

The methods of this paper were initially developed with a view towards bounding the regularity of subspace arrangements.  Recently, a sharp bound for the regularity of arbitrary subspace arrangements has been proved jointly by the author and Harm Derksen using more algebraic methods \cite{derksen-sidman}.  However, we shall see that in certain low-dimensional settings algebraic complexity is invisible from the perspective of sheaves, and this leads naturally to statements about the regularity of products and tensor products.

In \S 1 we define regularity and recall basic facts.  We then prove Lemma \ref{lem: homology} which is a technical statement regarding the regularity of the zero-th homology of a complex with low-dimensional homology that strengthens Lemma 1.6 in \cite{gruson-lazarsfeld-peskine}.   From this we get a more general bound on the regularity of tensor products than was previously known.  The applications to products of ideal sheaves and unions of schemes that intersect pairwise in points are part of the folklore of the subject and follow easily.  The main result of \S 1 is Theorem \ref{thm: prod} which tells us that these methods also give bounds on the regularity of the products of the homogeneous ideals when their intersection has dimension zero.  In \S 2 we give applications to the regularity of arrangements of two-planes based on the idea of making approximations with sheaves.

I am deeply indebted to Rob Lazarsfeld for introducing me to this area and for his guidance in my work.  I am very grateful to Bernd Sturmfels for his encouragement and to Aldo Conca and J\"urgen Herzog for their kind communications and for sending early versions of their preprints.  I also thank Tamon Stephen for help searching the combinatorial literature, Dan Rogalski and Tom Weston for helpful conversations, and the referee, Harm Derksen, and Mihnea Popa for comments on the exposition.  I thank NATO for their support of the \emph{Workshop in Exterior Algebras and other new directions in Algebraic Geometry, Commutative Algebra, and Combinatorics} held in Sicily in September of 2001.

\section{On regularity and products}

In this section we will show that regularity is subadditive for tensor products of coherent sheaves that are locally free away from low-dimensional sets and that the regularity of products of ideals and ideal sheaves is subadditive when the corresponding schemes have low-dimensional intersections.  We also include applications to subspace arrangements with singularities of dimension less than one.  The proofs exploit the fact that the support of a coherent sheaf is closed and that the higher cohomology groups vanish above the dimension of the support.  In other words, if \sheaf{F} is a coherent sheaf on \p{n} and  $\dim \supp$ $\sheaf{F} = m$, then $\coh{i}{\p{n}}{\sheaf{F}(l)} = 0$ for $i > m$ and all $l$.

We begin by recalling the definition of regularity for a finitely generated graded module over a polynomial ring:

\begin{defin}  Let $M$ be a finitely generated graded module over $S$ and let \[ 0 \to E_{n+1} \to \cdots \to E_1 \to E_0 \to M \to 0\] be a minimal free resolution of $M$.  Then $M$ is $m$-regular if $E_i$ is generated in degrees less than or equal to $m+i$ and the regularity of $M$, denoted $\reg(M)$, is the least $m$ for which this holds.  We will say that a projective scheme $X$ is $m$-regular if its saturated homogeneous ideal is $m$-regular and that the regularity of $X$ is the regularity of this ideal.
\end{defin}

There is a corresponding notion of regularity for coherent sheaves on \p{n}:

\begin{defin}  Let $\sheaf{F}$ be a coherent sheaf on \p{n}.  Then \sheaf{F} is $m$-regular if \[\coh{i}{\p{n}}{\sheaf{F}(m-i)} = 0\] for all $i > 0,$ and we say that the regularity of \sheaf{F} is the least $m$ for which these vanishings hold.
\end{defin}

The relationship between the regularity of an ideal and the regularity of the corresponding coherent sheaf of ideals is given by (Definition 3.2 in \cite{bayer-mumford}, with proof in their technical appendix):

\begin{thm}Let $I$ be an ideal of $S$ and \sheaf{I} be the corresponding sheaf.  Then the following properties are equivalent:
\label{thm: regdef}
\begin{itemize}

\item[(a)] the natural map $I_m \to \coh{0}{\p{n}}{\sheaf{I}(m)}$ is an isomorphism and \coh{i}{\p{n}}{\sheaf{I}(m-i)} = 0, $1 \leq i \leq n.$

\item[(b)] the natural maps  $I_d \to \coh{0}{\p{n}}{\sheaf{I}(d)}$ are isomorphisms for all $d \geq m$ and  \coh{i}{\p{n}}{\sheaf{I}(d)} = 0, $d + i \geq m$, $i \geq 1.$

\item[(c)] Take a minimal resolution of $I$ by free graded $S$-modules:
\[0 \to \oplus^{r_n}_{\alpha = 1}S(-e_{\alpha, n}) \to \cdots \to \oplus{r_0}_{\alpha = 1} S(- e_{\alpha, 0}) \to I \to 0\] Then $deg(e_{\alpha, i}) \leq m+i $ for all $\alpha, i$.  
\end{itemize} 

\end{thm}

Note that if a finitely generated graded module $M$ is $m$-regular then the associated sheaf is also $m$-regular and that any module agreeing with $M$ in degrees $k$ and higher for $k \ge m$ is $k$-regular.

One should also remark that it is well known (see \cite{eisenbud-goto}) that if $M$ is $m$-regular then the truncated module $M_{\ge m}$ has an $m$-linear resolution \[ \cdots \to \oplus S(-m-2) \to \oplus S(-m-1) \to \oplus S(-m) \to M_{\ge m} \to 0. \]  This implies that the sheaf $\sheaf{M}$ associated to $M$ has an $m$-linear resolution \[ \cdots \to \oplus \struct{\p{n}}(-m-2) \to \oplus \struct{\p{n}}(-m-1) \to \oplus  \struct{\p{n}}(-m) \to \sheaf{M} \to 0.\]

The following lemma, which is a variant of Lemma 1.6  in \cite{gruson-lazarsfeld-peskine}, is key in allowing us to control the regularity of tensor products of sheaves with higher $\tor$'s supported on sets of low dimension.

\begin{lem}
\label{lem: homology}
Let 
\[
\xymatrix{
\dots \ar[r]^{\phi_3} & \sheaf{E}_2 \ar[r]^{\phi_2} & \sheaf{E}_1 \ar[r]^{\phi_1} &\sheaf{E}_0 \ar[r]& 0
}
\]
\noindent be a complex of sheaves on \p{n} with homology sheaves $\sheaf{H}_i$, for $i\ge 0$.  Suppose that $ \sheaf{E}_i$ is $(m+i)$-regular and the dimension of the support of the higher homology of the complex is less than or equal to two, i.e., \[\dim \supp \sheaf{H}_i \leq 2 \mbox{ for } i\ge 1.\]  Then $\sheaf{H}_0 := \im(\phi_1)$ is $m$-regular.
\end{lem}

\begin{proof}  We need to show that \coh{i}{\p{n}}{\sheaf{H}_0(m-i)} vanishes for all $i > 0$.  Consider the following diagram: \[\xymatrix{
&& 0& 0 \ar[d]\\
&0 \ar[r]&\sheaf{B}_1 \ar[r] \ar[u] & \sheaf{Z}_1 \ar[d] \ar[r]& \sheaf{H}_1 \ar[r] &0\\
&\cdots \ar[r] & \sheaf{E}_2 \ar[r] \ar[u]&\sheaf{E}_1 \ar[r]\ar[d] & \sheaf{E}_0 \ar[r]  & 0\\
&& 0 \ar[r] &  \sheaf{B}_0 \ar[r] \ar[d] & \sheaf{Z}_0 \ar[r] \ar[u]& \sheaf{H}_0 \ar[r] & 0\\
& &&0 & 0 \ar[u] 
}\] where $\sheaf{Z}_i := \ker(\phi_i)$ and $\sheaf{B}_i := \im(\phi_{i+1})$.

Fix $i \ge 1$ and twist the entire diagram by $(m-i)$.  Examining the long exact sequence \[\cdots \to \coh{i}{\p{n}}{\sheaf{Z}_0(m-i)} \to \coh{i}{\p{n}}{\sheaf{H}_{0}(m-i)} \to \coh{i+1}{\p{n}}{\sheaf{B}_0(m-i)} \to \cdots,\] we see that we need to show that $\coh{i}{\p{n}}{\sheaf{Z}_0(m-i)} = 0$ and $\coh{i+1}{\p{n}}{\sheaf{B}_0(m-i)} = 0$.

Run the long exact sequence associated to the vertical and horizontal short exact sequences in the diagram.  Note that  $\coh{i+j}{\p{n}}{\sheaf{E}_j(m-i)} = 0$ for all $j \ge 0$ because $\sheaf{E}_j$ is $(m+j)$-regular, and $\coh{i+j+1}{\p{n}}{\sheaf{H}_j(m-i)} = 0$ for all $j \ge 1$ because $\sheaf{H}_j$ is supported on a set of dimension $\leq 2$.  Additionally, $\coh{i+j+1}{\p{n}}{\sheaf{Z}_j(m-i)} = 0$ for $j = n-i$ from dimension considerations.  These vanishings give the vanishing of each $\coh{i+j+1}{\p{n}}{\sheaf{Z}_j(m-i)}$ and $\coh{i+j}{\p{n}}{\sheaf{B}_{j-1}(m-i)}$ for $j \ge 1$.  Finally, we know that $\coh{i}{\p{n}}{\sheaf{Z}_0(m-i)} = 0$ because $\sheaf{Z}_0 = \sheaf{E}_0$ is $m$-regular by hypothesis.  Therefore, we conclude that $\coh{i}{\p{n}}{\sheaf{H}_{0}(m-i)}$ vanishes and that $\sheaf{H}_0$ is $m$-regular.
\end{proof}

The primary application of Lemma \ref{lem: homology} for our purposes is the following:

\begin{prop}
Let  \sheaf{F} and \sheaf{G} be coherent sheaves of \struct{\p{n}}-modules that are $f$ and $g$-regular, respectively.  If there is a subvariety $V \subset \p{n}$ of dimension less than or equal to two such that at least one of \sheaf{F} or \sheaf{G} is locally free at any point not on $V$, then $\sheaf{F}\otimes\sheaf{G}$ is $(f+g)$-regular.
\label{prop: tensorreg}
\end{prop}

\begin{proof} Recall that $\tor_*(\sheaf{F}, \sheaf{G})$ can be computed by resolving either \sheaf{F} or \sheaf{G}.  Since \sheaf{F} is $f$-regular we can take an $f$-linear resolution of vector bundles:
\[
\xymatrix{
 & \sheaf{E}_2  \ar@{=}[d] & \sheaf{E}_1 \ar@{=}[d] & \sheaf{E}_0 \ar@{=}[d] \\
\cdots \ar[r] & \oplus \struct{\p{n}}(-f-2) \ar[r] & \oplus \struct{\p{n}}(-f-1) \ar[r] & \oplus \struct{\p{n}}(-f) \ar[r] &  \sheaf{F} \ar[r] &0.
}
\]
\noindent Note that $\sheaf{H}_i(\sheaf{E}. \otimes \sheaf{G}) = \tor_i(\sheaf{F}, \sheaf{G}).$  Over an open set $U$ where \sheaf{G} is free, tensoring by \sheaf{G} is exact, which implies that the higher homology of $\sheaf{E}.\otimes \sheaf{G}$ vanishes on $U$.  By reversing the roles of $\sheaf{F}$ and $\sheaf{G}$ and applying the same argument we see that these higher homology sheaves also vanish locally wherever \sheaf{F} is free.  Therefore, by our hypotheses the higher homology of the complex $\sheaf{E}. \otimes \sheaf{G}$ must be supported only on $V$.  The result now follows from Lemma \ref{lem: homology} since each $\sheaf{E}_i \otimes \sheaf{G}$ is $(f+g+i)$-regular and $\sheaf{H}_0(\sheaf{E}.\otimes \sheaf{G})$ is $\sheaf{F}\otimes\sheaf{G}$.
\end{proof}

As a corollary to Proposition \ref{prop: tensorreg} we have the following statement about the regularity of the union of disjoint schemes:

\begin{cor}  Let $X_1, \ldots, X_d$ be $d$ pairwise disjoint schemes in \p{n} and let $m_i$ denote the regularity of $X_i$.  Then the ideal sheaf of their union is $\sum m_i$-regular.
\label{cor: disjoint}
\end{cor} 

\begin{proof}  Let $\sheaf{I}_1, \ldots, \sheaf{I}_d$ be the ideal sheaves of the $d$ schemes.  Using Proposition \ref{prop: tensorreg} we can see by induction that $\sheaf{I}_1 \otimes \cdots \otimes \sheaf{I}_d$ is $\sum m_i$-regular.  To see that $\sheaf{I}_1 \otimes \cdots \otimes \sheaf{I}_d$ is equal to $\sheaf{I}(\cup X_i) = \sheaf{I}_1 \cap \cdots \cap \sheaf{I}_d$ we examine the stalks of the sheaves.  At any closed point $p$ at most one of the sheaves, without loss of generality, $\sheaf{I}_1$, is nontrivial.  Therefore, \[(\sheaf{I}_1 \otimes \cdots \otimes \sheaf{I}_d)_p = (\sheaf{I}_1)_p\otimes (\struct{\p{n}})_p \otimes \cdots \otimes (\struct{\p{n}})_p = (\sheaf{I}_1)_p\] which is clearly equal to \[(\sheaf{I}_1 \cap \cdots \cap \sheaf{I}_d)_p = (\sheaf{I}_1)_p \cap (\struct{\p{n}})_p \cap \cdots \cap (\struct{\p{n}})_p = (\sheaf{I}_1)_p.\]
\end{proof}

\rmk  Corollary \ref{cor: disjoint} shows that the ideal sheaf of a scheme $X$ consisting of $d$ disjoint linear spaces in \p{n} is $d$-regular.  (For an arrangement of disjoint subspaces one can also argue directly that the $(d-1)$-forms on \p{n} surject onto $\Gamma(X,\struct{X}(d-1))$ to achieve the vanishing of $\coh{1}{\p{n}}{\sheaf{I}(d-1)},$ which is the crucial vanishing.)  In particular, we have that the ideal sheaf of $d$ distinct points in \p{n} is $d$-regular.

\bigskip

Recall that the tensor product of ideal sheaves is not an ideal sheaf itself.  However, it maps surjectively onto the product of the ideal sheaves, which is again an ideal sheaf.  The following lemma shows that under suitable hypotheses on dimension, the regularity of the product of ideal sheaves is subadditive.

\begin{prop} 
\label{prop: productreg}
Suppose that \sheaf{I} and \sheaf{J} are ideal sheaves on \p{n} with regularities $m_1$ and $m_2$, respectively.  Suppose also that the zeros of \sheaf{I} and \sheaf{J} intersect in a set of dimension at most one.  Then $\sheaf{I} \cdot \sheaf{J}$ is $(m_1+m_2)$-regular.
\end{prop}

\begin{proof}  First, notice that the hypotheses of Proposition \ref{prop: tensorreg} above are satisfied.  So $\sheaf{I} \otimes \sheaf{J}$ is $(m_1+m_2)$-regular.  Now, since $\sheaf{I} \cdot \sheaf{J}$ is the image of $\sheaf{I} \otimes \sheaf{J}$ in \struct{X}, we have the following short exact sequence of sheaves \[ \ses{\sheaf{K}}{\sheaf{I} \otimes \sheaf{J}}{\sheaf{I} \cdot \sheaf{J}},\] and the support of \sheaf{K} is contained in a set of dimension at most one.  Running the long exact sequence, we see that $\coh{i}{\p{n}}{\sheaf{I} \cdot \sheaf{J}(m_1+m_2 - i)} = 0$ for all $i > 0$ since $\coh{i}{\p{n}}{\sheaf{K}(l)}$ vanishes for $ i > 1$ and any $l.$ 
\end{proof}

In general, Propostion \ref{prop: productreg} does not say anything about the regularity of products of ideals.  However, when the zero sets of the ideals intersect in a finite set of points, we do have the corresponding result for the product of the ideals themselves.

\begin{thm} \label{thm: prod} 
If $I$ and $J$ are homogeneous ideals of $S$ defining schemes in \p{n} whose intersection is a finite set of points, then $\reg(I \cdot J) \leq \reg(I)+\reg(J)$.
\end{thm}
\begin{proof}
Let $m_1 = \reg(I)$ and $m_2 = \reg(J)$.  Let $\sheaf{I}$ be the sheaf associated to $I$ and \sheaf{J} be the sheaf associated to $J$.  By Proposition \ref{prop: productreg} and part (b) of Theorem \ref{thm: regdef} it is enough to show that $I \cdot J$ agrees with its saturation in degrees greater than or equal to $m_1+m_2$.  Since $I \cdot J$ is an ideal, this is equivalent to showing that $(I \cdot J)_d$ surjects onto $\coh{0}{\p{n}}{\sheaf{I}\cdot \sheaf{J}(d)}$ for $d \ge m_1+m_2$.  

For $m \ge m_1$ let \sheaf{P}. be an $m$-linear resolution of \sheaf{I} with \[\sheaf{P}_0 := I_m \otimes_k \struct{\p{n}}(-m),\] and tensor the resolution by \sheaf{J}.  The essential point is that because the higher homology sheaves of the resulting complex have zero-dimensional support, the map \[ \coh{0}{\p{n}}{\sheaf{P}_0 \otimes \sheaf{J}(m+m_2)} \to \coh{0}{\p{n}}{\sheaf{I}\otimes\sheaf{J}(m+m_2)} \] is a surjection.  But our definition of $\sheaf{P}_0$ implies that $ \coh{0}{\p{n}}{\sheaf{P}_0 \otimes \sheaf{J}(m+m_2)}$ equals $I_{m} \otimes_k J_{m_2}$ so we have the surjection \[I_{m} \otimes_k J_{m_2} \to \coh{0}{\p{n}}{\sheaf{I}\otimes\sheaf{J}(m+m_2)}.\]  Additionally, our hypotheses imply that  $\coh{0}{\p{n}}{\sheaf{I}\otimes\sheaf{J}(l)}$ surjects onto $\coh{0}{\p{n}}{\sheaf{I}\cdot\sheaf{J}(l)}$ for any $l$ because the support of the kernel of the map from $\sheaf{I} \otimes \sheaf{J}$ to $\sheaf{I} \cdot \sheaf{J}$ is zero-dimensional.  Composing the maps we see that  $I_{m} \otimes_k J_{m_2} \to \coh{0}{\p{n}}{\sheaf{I} \cdot \sheaf{J}(m+m_2)}$ is surjective, and since its image is $I_{m} \cdot J_{m_2} = (I \cdot J)_{m+m_2}$, we are done.  
\end{proof}
Using Theorem \ref{thm: prod} we can recover the bound on the regularity of powers of an ideal of a finite set of points given by Chandler \cite{chandler} and independently by Geramita, Gimigliano, and Pitteloud \cite{geramita-gimigliano-pitteloud} by taking $I = J$ and working inductively.  Additionally, Conca and Herzog have recently proved a similar bound for the regularity of the product of an ideal $I$ with any module if $\dim(S/I) \leq 1$  \cite{conca-herzog}, which we can recover in the special case where the module is an ideal.   Note that bounds of this type fail to hold in general.  Sturmfels \cite{sturmfels} and Terai (Remark 3 in \cite{conca}) have given examples of ideals $I$ for which $\reg(I^2) > 2 \cdot \reg(I)$, and Conca and Herzog have shown that the regularity of a product of distinct ideals may be strictly greater than the sum of the regularities of its factors in \cite{conca-herzog}.

As a consequence of Proposition \ref{prop: productreg} we have a statement about the regularity of the union of schemes intersecting in points:

\begin{cor}
Let \sheaf{I} be the ideal sheaf of a projective scheme $X$ that consists of the union of $d$ schemes $X_1, \ldots, X_d$ in \p{n} whose pairwise intersections are finite sets of points.  Let $m_i$ be the regularity of $X_i$.  Then \sheaf{I} is $\sum m_i$-regular.
\end{cor}

\begin{proof} Let $\sheaf{I}_i$ be the ideal sheaf of $X_i$.  Consider the short exact sequence \[ 0 \to \sheaf{I}_1 \cdots \sheaf{I}_{d} \to \sheaf{I} \to \sheaf{C} \to 0.\]  The sheaf $ \sheaf{I}_1 \cdots \sheaf{I}_{d}$ is $\sum m_i$-regular by an induction argument applied to Proposition \ref{prop: productreg} and it differs from \sheaf{I} only where the schemes intersect.  Since all of the intersections are zero-dimensional by hypothesis the support of $\sheaf{C}$ must be zero-dimensional.  But then the higher cohomology of any twist of \sheaf{C} is automatically zero, and it follows that \sheaf{I} is $\sum m_i$-regular
\end{proof}

This give us the special case:

\begin{cor}
The ideal sheaf of $d$ $k$-planes in \p{n} meeting only in points is $d$-regular.  (In particular the union of $d$ lines is $d$-regular.)
\end{cor}

\section{Cones and 2-planes}
This section is based on the idea that if a morphism of sheaves $\sheaf{F} \to \sheaf{G}$ produces a kernel and cokernel with support of lower dimension than the support of the original sheaves, then we can view \sheaf{F} as an approximation of \sheaf{G} and attempt to use this approximation to study the regularity of \sheaf{G}.  We will use this idea to show that the ideal sheaf of a collection of $d$ linear spaces in \p{n} with pairwise intersections of dimension less than or equal to one, is $d$-regular.  Our starting point is that the ideal sheaf of a subspace arrangement is equal to the product of the ideal sheaves of its constituent subspaces except where they intersect.  We use cones over lower dimensional subspace arrangements to model these intersections.

We begin by fixing notation.  Let $X$ be the union of $d$ linear spaces $X_1, \ldots, X_d$ with ideal sheaves $\sheaf{I}_1, \ldots, \sheaf{I}_d$, such that the pairwise intersections have dimension at most one.  Also let \sheaf{I} be the ideal sheaf of $X$ and $\sheaf{J} := \sheaf{I}_1 \cdots \sheaf{I}_d$.  

Our first approximation comes from \[ 0 \to \sheaf{J} \to \sheaf{I} \to \frac{\sheaf{I}}{\sheaf{J}} \to 0.\]  By Proposition \ref{prop: productreg} we know that \sheaf{J} is $d$-regular  because it is the product of $d$ ideals, each of which is 1-regular, and the planes corresponding to these ideals have pairwise intersections of dimension at most one.  More generally, Conca and Herzog have recently shown that the product of any $d$ linear ideals is $d$-regular \cite{conca-herzog}, and the $d$-regularity of \sheaf{J} follows by considering the sheafification of the product.  Comparing the stalks of $\sheaf{J}$ and $\sheaf{I}$ we see they agree everywhere except at the intersections of the planes, so the support of $\frac{\sheaf{I}}{\sheaf{J}}$ consists of lines and isolated points.  Heuristically, \sheaf{J} approximates \sheaf{I} because our map is an isomorphism over a dense open set.  The result will follow if we can show that $\frac{\sheaf{I}}{\sheaf{J}}$ is $d$-regular.  We will do this by exhibiting a $d$-regular sheaf that approximates it away from a zero-dimensional set.

Now we construct the sheaf that will approximate $\frac{\sheaf{I}}{\sheaf{J}}$. We will see that it is enough to approximate $\frac{\sheaf{I}}{\sheaf{J}}$ over the lines of its support.   Let $L$ be a line contained in the support of $\frac{\sheaf{I}}{\sheaf{J}}.$ Define 
\[ \sheaf{I}_L := \sheaf{I}(\underset{X_i \supseteq L}{\cup} X_i), \hspace{.25in}  \sheaf{J}_L := \prod_{X_i \supseteq L} \sheaf{I}_i,\] \[ \sheaf{J}_L' := \prod_{X_i \nsupseteq L} \sheaf{I}_i, \hspace{.1in} \mbox{ and } \hspace{.1in}  d_L := \#\{X_i : X_i \supseteq L\}.\]

Therefore, there are $d_L$ planes passing through $L$ and $d-d_L$ planes that have either empty or zero-dimensional intersections with $L$.  We know that $\sheaf{J}_L'$ is $(d-d_L)$-regular by the same argument used above for \sheaf{J}.  Furthermore, in this case, $\sheaf{I}_L$ is $d_L$-regular, which follows from the lemma below whose proof we leave to the reader.

\begin{lem}\label{lem: cone-reg} If $\sheaf{I}$ is the ideal sheaf of a scheme $V$ in \p{n} with regularity $d$ and $\sheaf{I}^+$ is the sheaf of the cone over $V$ in \p{{n+m}}, then $\sheaf{I}^+$ is $d$-regular.
\end{lem}

To finish our construction, notice that the zero locus of $\sheaf{J}_L'$ intersects the zero locus of $\sheaf{I}_L$ in a set of dimension at most one.  Therefore, Proposition \ref{prop: productreg} applies and $\sheaf{I}_L \cdot \sheaf{J}_L'$ is $d$-regular.  Knowing this, we see that $\cok$ is $d$-regular.  In Proposition \ref{prop: directsum} we will show that if we sum over all lines of the support of $\frac{\sheaf{I}}{\sheaf{J}}$, $\oplus \cok$ approximates $\frac{\sheaf{I}}{\sheaf{J}}$ away from points.  For this we will need Lemma \ref{lem: support} which shows that each $\cok$ is supported only on $L$. 

\begin{lem}
\label{lem: support}
For any line $L$ in the support of $\frac{\sheaf{I}}{\sheaf{J}}$ the support of  $\cok$ is contained in $L$.
\end{lem}

\begin{proof}
Fix a point $q$ not on $L$.  We will show that the stalk of\[ \cok = \frac{\sheaf{I}_L \cdot \sheaf{J}_L'}{\sheaf{J}_L \cdot \sheaf{J}_L'}\] is zero at $q$.  The main point is that at $q$ the stalk of $\sheaf{I}_L$ is equal to the stalk of $\sheaf{J}_L$.  This follows from the observation that if $q$ is a point not on $L$ then at most one of the planes containing $L$ can pass through $q$.  There are two cases.  In the first case some plane $X_i$ containing $L$ also passes through $q$.  Then both the stalks of $\sheaf{I}_L$ and $\sheaf{J}_L$ at $q$ are equal to the stalk of $\sheaf{I}_i$ at $q$.  In the other case $q$ is not on any of the planes that pass through $L$.  Then the stalks of $\sheaf{I}_L$ and $\sheaf{J}_L$ are both trivial at $q$.  
\end{proof}

Notice that the induced vertical map on the right in the diagram below is an inclusion.  
\[
\xymatrix{
 0 \ar[r] & \sheaf{J} \ar[r] \ar[d] & \sheaf{I}_L \cdot \sheaf{J}_L' \ar[r] \ar[d] & \displaystyle{\cok} \ar[r]\ar@{-->}[d]  & 0\\
 0 \ar[r] & \sheaf{J} \ar[r] & \sheaf{I}  \ar[r] & \displaystyle{\frac{\sheaf{I}}{\sheaf{J}}} \ar[r] & 0
}
\]
\bigskip

Proposition \ref{prop: directsum} shows that $\oplus \cok$ is indeed an approximation of $\frac{\sheaf{I}}{\sheaf{J}}$ away from a zero-dimensional set.

\begin{prop}
\label{prop: directsum}
In the exact sequence  \[ 0 \to \sheaf{K} \to \oplus \cok \to \frac{\sheaf{I}}{\sheaf{J}} \to \sheaf{S} \to 0\]  \sheaf{K} and \sheaf{S} are skyscraper sheaves.
\end{prop}

\begin{proof}
Fix a line $L$ in the direct sum, and define $U_L$ to be \p{n} minus all of the $X_i$ that do not contain $L$.  When we restrict to $U_L$ all of the relevant ideal sheaves behave as if only the planes passing through $L$ exist.  More specifically, $\sheaf{J}_L'|_{U_L}$ is trivial and $\sheaf{I}|_{U_L}$ is equal to $\sheaf{I}_L|_{U_L}$.  We conclude that $\cok$ restricts to $\frac{\sheaf{I}}{\sheaf{J}}$ over $U_L$.  

To see that \sheaf{S} is a skyscraper sheaf, note that $U_L \cap L$ is a dense open subset of $L$.  Therefore, the intersection of the support of \sheaf{S} with $L$ is a finite set of points. Since this is true for each $L$, the support of \sheaf{S} is zero-dimensional.

Finally, because each of the terms in the direct sum $\oplus \cok$ includes into $\frac{\sheaf{I}}{\sheaf{J}}$, the kernel of the map from their direct sum may be supported only where their supports intersect.  By Lemma \ref{lem: support} each $\cok$ is supported only on $L$ so the intersections of the supports of the sheaves in the direct summand must be isolated points.  Therefore, the support of \sheaf{K} is also zero-dimensional.
\end{proof}

We now prove the main result of this section:

\begin{thm}
\label{thm: one-intersect}  If $\sheaf{I}$ is the ideal sheaf of $d$ linear spaces in \p{n} with pairwise intersections of dimension less than or equal to one, then \sheaf{I} is $d$-regular.
\end{thm}

\begin{proof}
Using Proposition \ref{prop: directsum} we have an exact sequence \[ 0 \to \sheaf{K} \to \oplus \cok \to \frac{\sheaf{I}}{\sheaf{J}} \to \sheaf{S} \to 0\] where \sheaf{K} and \sheaf{S} are skyscraper sheaves.  We already know that all of the sheaves are $d$-regular except for $\frac{\sheaf{I}}{\sheaf{J}}$.  But this follows from splitting the sequence into two short exact sequences and running the corresponding long exact sequences.  The result then follows by running the long exact sequence corresponding to \[\ses{\sheaf{J}}{\sheaf{I}}{\frac{\sheaf{I}}{\sheaf{J}}}. \]
\end{proof}

As a corollary to Theorem \ref{thm: one-intersect} we have the case of 2-planes:

\begin{cor}
\label{cor: 2-planes}
If $\sheaf{I}$ is the ideal sheaf of $d$ 2-planes in \p{n} then $\sheaf{I}$ is $d$-regular.
\end{cor}

\bibliographystyle{amsplain}

\end{document}